\documentclass[11pt,english]{smfart}
\usepackage[T2A]{fontenc}
\usepackage[utf8]{inputenc}
\usepackage[english,russian]{babel}

\usepackage{url}
\usepackage{smfthm}
\usepackage{hyperref}
\usepackage{MnSymbol}

\usepackage{graphicx}
\usepackage[matrix, arrow,all,cmtip,color]{xy}

\newcommand{\bi}{\begin{itemize}}
\newcommand{\ei}{\end{itemize}}
\newcommand{\bd}{\begin{description}}
\newcommand{\ed}{\end{description}}

\def\lra{\longrightarrow}

\def\rtt{\,{\Huge{\rightthreetimes}}\,}
\def\rttt{\,{\underline{\rightthreetimes}}\,}
\def\xra{\xrightarrow}

\def\Ker{{\text{Ker}\,}}
\def\Imm{{\text{Im}\,}}

\def\ZZ{\Bbb Z}

\def\nto{\not\to}
\def\ZpZ{{\Bbb Z}/\!p{\Bbb Z}}

\def\Ab{\textit{Ab}}\def\AbKer{{\textit{AbKer}}}

\def\lr{{\rtt lr}} 
\def\rr{{\rtt rr}}

\begin{document}

\title[
Basic finite group theory via the lifting property
]{
Formulating basic notions of \\ finite group theory\\ via the lifting property}
\author[masha gavrilovich]{masha gavrilovich
\thanks{\tiny Institute for Regional Economic Studies, Russian Academy of Sciences (IRES RAS). 
National Research University Higher School of Economics, Saint-Petersburg.
{\tt \tiny mi\!\!\!ishap\!\!\!p@sd\!\!\!df.org \url{http://mishap.sdf.org}}.\\ 
This paper commemorates the centennial of the birth of N.A.~Shanin, 
the teacher of S.Yu.Maslov and G.E.Mints, who was my teacher. I hope the motivation behind this paper is in spirit of the Shanin's group ТРЭПЛО (теоритическая разработка эвристического поиска логических обоснований, theoretical development of heuristic search for logical evidence/arguments).}
}
\begin{abstract}
We reformulate several basic notions of notions in finite group theory in terms
of iterations of the lifting property (orthogonality) with respect to
particular morphisms. Our examples include the notions 
being nilpotent,
solvable, perfect, torsion-free; p-groups and prime-to-p-groups; Fitting
subgroup, perfect core, p-core, and prime-to-p core. 
 
We also reformulate as in similar terms the conjecture that a localisation of a
(transfinitely) nilpotent group is (transfinitely) nilpotent.

\end{abstract}
\maketitle
\section{Introduction.}

We observe that several 
standard
elementary notions of finite group theory 
 can be defined by 
iteratively applying the same diagram chasing ``trick'',
namely the lifting property (orthogonality of morphisms), 
to simple classes of homomorphisms of finite groups.
%

The notions include a finite group being nilpotent, solvable, perfect, 
torsion-free; $p$-groups, and prime-to-$p$ groups; $p$-core, the Fitting subgroup, cf.\S\ref{sec:neg}-\ref{sec:reforms}. 

In \S\ref{sec:locals} we reformulate as a labelled commutative diagram the conjecture 
that a localisation
of a transfinitely nilpotent group is transfinitely nilpotent; this suggests a variety of related questions
and is inspired by the conjecture of Farjoun that a localisation of a nilpotent group is nilpotent.

The goal of this paper to present a collection of examples which show 
the lifting property is all that's needed to be able to
define a number of notions from simplest (counter)examples of interest.

Curiously, our observations lead to a concise and uniform notation (Theorem~\ref{theo:1}, Corollary~\ref{coro:1}~and~\ref{coro:2}), e.g.
 $$   (\ZpZ\lra 0)^{rr},\ \ (\AbKer)^{lr},\ \text{ and }\  (0\lra *)^{lr} $$ 
denote the classes of homomorphisms (of finite groups) 
 whose kernel is a $p$-group, resp.~soluble, subgroup, and those corresponding to subnormal subgroups.
One might hope that a notation so concise and uniform might be of use in computer algebra and automated theorem provers.

Deciphering this notation can be used as an elementary exercise in a first course of group theory or category theory on basic definitions and diagram chasing.




Such reformulations lead one to the following questions:

\bi\item
Can one extend this notation to capture more of finite group theory? 
\item Is this a hint towards 
category theoretic point of view on finite group theory? 
\ei

If one believes the evidence
provided by our examples 
is strong enough to demand an explanation, then
one should perhaps start by trying to find more examples defined in this way, and 
by calculating the classes of homomorphisms obtained by 
iteratively applying the Quillen lifting
property to simple classes of morphisms of finite groups.

\subsubsection*{Motivation} Our motivation was to formulate 
part of finite group theory in a form amenable 
to automated theorem proving while remaining human readable;
[GP] tried to do the same thing for the basics of general topology. 



Little attempt has been made to go beyond these examples. 
Hence open questions remain: are there other interesting examples of lifting properties in the category of (finite) groups?
Can a complete group-theoretic argument be reformulated in terms of diagram chasing, say 
the classification of CA-groups or $pq$-groups, or 
elementary properties of subgroup series; can category theory notation
be used to make expositions easier to read? 
Can these reformulations be used in automatic theorem proving? 
Is there a decidable fragment of (finite) group theory
 based on the Quillen lifting property and, more generally,
 diagram chasing, cf.~[GLS]?
Can the Sylow theorems (only existence and uniqueness of Sylow subgroups)
be proven using this characterization of p-groups?
Could the components of a finite group, and their properties (commute
pairwise, commute with normal p-subgroups) be characterized and proven to exist with
these methods?

\section{Definitions and examples of reformulations} 

\subsection{Key definition: the Quillen lifting property (negation/orthogonal)} 
The Quillen lifting property, also known as orthogonality of morphisms, is a property of a pair of morphisms in a category. 
It appears in a prominent way
   in the theory of model categories, an axiomatic framework for homotopy theory introduced by Daniel Quillen,
and 
is used to define properties of morphisms starting from an explicitly given class of morphisms. 

%

\begin{defi} We say that two morphisms $A \xra f B$ and $X \xra g Y$ in a category $C$ are {\em orthogonal} and write $f\rtt g$
iff 
for each $i:A\lra X$ and $j:B\lra Y$
making the square commutative, i.e.~$f\circ j=i\circ g$ there is a diagonal arrow $\tilde j:B\lra X$ making the total diagram
$A\xra f B\xra {\tilde j} X\xra g Y, A\xra i X, B\xra j Y$ commutative, i.e.~$f\circ \tilde j=i$ and $\tilde j\circ g=j$ (see Figure 1a).

We may also say that {\em $f$ lifts wrt $g$}, {\em $f$ left-lifts wrt $g$}, or {\em $g$ right-lifts wrt $f$}, or that $f$ {\em antagonizes} $g$.

%

By analogy with orthogonal complement of a non-symmetric bilinear form, 
define  {\em left/right Quillen negation} or {\em left/right Quillen orthogonal} of a class $P$
of morphisms: 
 $$ P^{\rtt l}:=\{\, f\ :\ f\,\rtt\,g\,\ \text{for each }g\in P\,\}$$
$$P^{\rtt r}:=\{\, g\ :\ f\, \rtt g\,\ \text{for each }g\in P\,\}$$
\end{defi}
\vskip 3pt

We have 
$$ P^{\rtt l}=P^{\rtt lrl},\  P^{\rtt r}=P^{\rtt rlr},\  P\subset P^{\rtt lr},\  P\subset P^{\rtt rl}$$
$$P \subset Q \text{ implies } Q^{\rtt l} \subset P^{\rtt l},\ Q^{\rtt r} \subset P^{\rtt r},\  P^{\rtt lr} \subset Q^{\rtt lr},\ P^{\rtt rl} \subset Q^{\rtt rl} $$
$$ P\cap P^{\rtt l} \subset (Isom),\ \  P\cap P^{\rtt r}\subset (Isom)$$
Under certain assumptions on the category and property $P$ 
Quillen small object argument shows that each morphism $G\lra H$ decomposes both as 
$$G \xra{(P)^{\rtt lr}}\cdot \xra{(P)^{\rtt r}} H \text{ and } 
G \xra{(P)^{\rtt r}}\cdot \xra{(P)^{\rtt lr}} H .$$
Using the Quillen lifting property is perhaps the simplest way to define a class of morphisms {\em without} a given property
in a manner useful in a category theoretic diagram chasing computation.

\def\rrt#1#2#3#4#5#6{\xymatrix{ {#1} \ar[r]^{} \ar@{->}[d]_{#2} & {#4} \ar[d]^{#5} \\ {#3}  \ar[r] \ar@{-->}[ur]^{}& {#6} }}
\begin{figure}
\begin{center}
\large
$(a)\  \xymatrix{ A \ar[r]^{i} \ar@{->}[d]_f & X \ar[d]^g \\ B \ar[r]_-{j} \ar@{-->}[ur]^{{\tilde j}}& Y }$
$\ \ \ \ \ \ (b)\ \xymatrix{ A \ar[r] \ar@{->}[d]_{(P)} & X \ar[d]^{\therefore (Q)} \\ B \ar[r] \ar@{-->}[ur]& Y }$
$\ \ \ \ \ \ (c)\ \xymatrix{ A \ar[r] \ar@{->}[d]_{\therefore(P)} & X \ar[d]^{ (Q)} \\ B \ar[r] \ar@{-->}[ur]& Y }$

\end{center}
\caption{\label{fig1}\small
(a) 
The definition of a lifting property $f\rtt g$: for each $i:A\lra X$ and $j:B\lra Y$
making the square commutative, i.e.~$f\circ j=i\circ g$, there is a diagonal arrow $\tilde j:B\lra X$ making the total diagram
$A\xra f B\xra {\tilde j} X\xra g Y, A\xra i X, B\xra j Y$ commutative, i.e.~$f\circ \tilde j=i$ and $\tilde j\circ g=j$.
We say that $f$ lifts wrt $g$, $f$ left-lifts wrt $g$, or $g$ right-lifts wrt $f$.
\ (b) Right Quillen negation. 
The diagram defines a property $Q$ of morphisms in terms of a property $P$; a morphism has property (label) $Q$ iff it right-lifts 
wrt any morphism with property $P$, i.e. 
$Q=\{\, q :\ p\, \rtt q\,\ \text{for each }p\in P\,\}$
\ (c) Left Quillen negation. 
The diagram defines a property $P$ of morphisms in terms of a property $Q$; a morphism has property (label) $P$ iff it left-lifts 
wrt any morphism with property $Q$, i.e. 
 $ P=\{\, p\ :\ p\,\rtt\,q\,\ \text{for each }q\in Q\,\}$
}
\end{figure}

\subsection{\label{sec:neg}A list of iterated Quillen negations of simple classes of morphisms}

Let $(0\lra *)$, resp.~$(0\lra Ab)$,  denote the class of morphisms from the trivial groups to an arbitrary group, resp. Abelian group.
Let $(*\lra 0)$, resp.~$(Ab\lra 0)$,  denote the class of morphisms to the trivial groups from an arbitrary group, resp. Abelian group.
Let $(AbKer)$ denote the class of homomorphisms with an Abelian kernel.

\begin{theo}\label{theo:1}
In the category of Finite Groups,
\begin{enumerate}

\item $(\AbKer)^\lr$ is the class of homomorphisms whose kernel is solvable
\item $(0\lra *)^\lr$ is the class of subnormal subgroups
\item  $(0\lra \Ab)^\lr=(0\lra \Ab)^\lr=\{[G,G]\lra G: G\text{ is arbitrary}\}^\lr$ is the class of subgroups $H<G$ such that there is a chain of subnormal subgroups 
$H=G_0 \vartriangleleft G_1  \vartriangleleft \ldots  \vartriangleleft G_n =G$ such that $G_{i+1}/G_{i}$ is Abelian, for $i=0,...,n-1$.
\item $(0 \lra S)^\lr$ is the class of subgroups $H<G$ such that there is a chain of subnormal subgroups
$H=G_0 \vartriangleleft G_1  \vartriangleleft \ldots  \vartriangleleft G_n =G$ such that $G_{i+1}/G_{i}$ embeds into $S$, for $i=0,...,n-1$. 

\item  $(\ZZ/p\ZZ\longrightarrow 0)^{r}$ is the class of homomorphisms whose kernel has no elements of order $p$
\item  $(\ZZ/p\ZZ\longrightarrow 0)^{\rr}$ is the class of surjective homomorphisms whose kernel is a $p$-group 
\end{enumerate}

In the category of Groups,
\begin{enumerate}

\item $(*\lra 0)^{\rtt l}$ is the class of retracts
\item  $(0\lra *)^{\rtt r}$ is the class of split homomorphisms

\item $(0\longrightarrow \ZZ)^{r}$ is the class of surjections
\item $(\ZZ\lra 0)^{r}$ is the class of injections
\item     a group \ensuremath{F} is free iff
 $0\longrightarrow F$ is in  $(0\longrightarrow \ZZ)^{rl}$ 

\item a group $A$ is Abelian iff $A\lra 0$ is in $( \Bbb F_2 \lra \ZZ\times\ZZ)^{\rtt r}$


\item group $G$ can be obtained from $H$ by adding commutation relations, i.e.~the kernel of $H\lra G$ is generated by commutators $[h_1,h_2]$, $h_1,h_2\in H$,
iff 
 $H\lra G$ is in $( \Bbb F_2 \lra \ZZ\times\ZZ)^{rl}$

\item subgroup $H$ of $G$ is a normal span of substitutions in words $w_1,..,w_i$ of the free group $\Bbb F_n$ iff $G \lra G/H$ is in $( \Bbb F_n \lra \Bbb F_n/<w_1,...,w_i>)^{rl}$

\item $(\AbKer)^{\rtt l}$ is the class of homomorphisms whose kernel is perfect
\end{enumerate}

\end{theo}

\begin{proof} The proof is a matter of deciphering notation. 

Proof of item 1. First note that $P\lra 0\in (\AbKer)^{\rtt l}$ for a perfect group, and $P\lra
0 \rtt H\lra G$ implies $\Ker(H\lra G)$ is soluble. This means that $(\AbKer)^{\rtt lr}$
is contained in the class of maps whose kernel is soluble. On the other hand, 
any such map is a composition of maps $H/[S_n,S_n]\lra H/S_n$, ..., $H/[S_0,S_0]\lra H/S_0$, and $\Imm H\lra G$ where $S_{n+1}=[S_n,S_n],S_0=\Ker(H\lra G)$
is the descending derived series.

Now let us prove item 2. 
By definition, $A\lra B$ is in  $(0\lra *)^{\rtt l}$ iff $A\lra B
\rtt 0\lra G $ for any group $G$. Take $G=B/A^B$ to be the quotient of $B$ by
the normal closure of $A$, and $B\lra G$ to be the quotient map. This shows
that if $G=B/A^B$ is non-trivial, then the lifting property fails. 
 On the other hand, it is easy to check the lifting property holds that in a commutative square,
the map to $G$ factors via $B/A^B$, hence the lifting property holds if $B/A^B$ is trivial. 

Let $C \vartriangleleft D$ be a normal subgroup. The lifting property $A\lra B
\rtt 0\lra D/C $ implies $A\lra B
\rtt C\lra D $. Orthogonals are necessarily closed under composition, hence this implies that 
if $C$ is a subnormal subgroup of $D$, i.e.~there exists a series
if $C \vartriangleleft D_n \vartriangleleft D_{n-1} \vartriangleleft ... \vartriangleleft D_1 \vartriangleleft D$,
then the lifting property holds and $C\lra D$ is in $(0\lra *)^{lr}$. 

Now assume that $C$ is not subnormal in $D$ and let $C<C'$ be a minimal subnormal subgroup of $D$ containing $C$. 
Then $C\lra C'$ is in $(0\lra *)^{\rtt l}$ and the lifting property $C\lra C' \rtt C\lra D$ fails, as required. 

Items 5 and 6 use Cauchy theorem that a prime $p$ divides the order of a group iff the group has an element of order $p$.

 \end{proof}

\subsection{\label{sec:reforms}Concise reformulations in terms of iterated Quillen negation}

We use the observations above to concisely reformulate several elementary notions in finite group theory.

\begin{coro}\label{coro:1} In the category of Finite Groups,  
\begin{enumerate}
\item a finite group $S$ is soluble iff either of the following equivalent conditions holds: 
\bi\item $S\lra 0$ is in $(\AbKer)^\lr$
\item $0\lra S$ is in $(0\lra \Ab)^\lr$
\ei
\item a finite group $H$ is nilpotent iff 
\begin{itemize}\item the diagonal map $H\lra H\times H$, $x\mapsto (x,x)$, is in  $(0\lra *)^\lr$ 
\end{itemize}
\item the Fitting subgroup $F$ of $G$ is the largest subgroup such that 
\bi \item the diagonal map $F\lra G\times G$, $x\mapsto (x,x)$, is in  $(0\lra *)^\lr$
\ei



\item a finite group $H$ is a $p$-group iff one of the following equivalent conditions hold:
 \begin{itemize}
  \item $H\lra 0$ is in  $(\ZZ/p\ZZ\lra 0)^\rr$
 \item $0\lra H$ is in $(0\lra \ZZ/p\ZZ)^\lr$
\end{itemize}
\end{enumerate} 

In the category of  Groups,  
\begin{enumerate}
\item a group $G$ is torsion-free iff  iff $0\longrightarrow G$ is in $\{ n\ZZ\longrightarrow \ZZ: n>0 \}^{\rtt r}$
\item a subgroup \ensuremath{H<G} contains torsion and is pure iff $H\longrightarrow G$ is in $\{ n\ZZ\longrightarrow \ZZ: n>0 \}^{\rtt r}$
\item $H$ is a verbal subgroup of $G$ generated by substitutions in words $w_1,..,w_i$ in the free group $\Bbb F_n$ iff $H$ fits into an exact sequence 
$$H \lra G \xra {(\Bbb F_n \lra\, \Bbb F_n/<w_1,..,w_i>)^{\rtt rl}} \cdot \xra{(\Bbb F_n \lra\, \Bbb F_n/<w_1,..,w_i>)^{\rtt r}} 0,$$ or, equivalently, is the kernel of the corresponding homomorphism

\item a finite group $S$ is transfinitely soluble,
 i.e.~there exists an ordinal $\alpha$ such that $G^\alpha=0$, where $G^{\beta+1}=[G,G^\beta]$, and $G^\beta=\cap_{\gamma<\beta} G^\gamma$ whenever  $\beta\neq \gamma+1$ for $\gamma<\beta$,
 iff 
\bi\item $S\lra 0$ is in $(\AbKer)^\lr$
\ei
\item a group $G$ is transfinitely nilpotent, i.e.~~there exists an ordinal $\alpha$ such that $G^\alpha=0$, where $G^{\beta+1}=[G,G^\beta]$, and $G^\beta=\cap_{\gamma<\beta} G^\gamma$ whenever $\beta\neq \gamma+1$ for $\gamma<\beta$,
iff
\begin{itemize}\item the diagonal map $H\lra H\times H$, $x\mapsto (x,x)$, is in  $(0\lra *)^\lr$ 
\end{itemize}
\end{enumerate}
\end{coro}

\begin{coro}
The statement that a group of odd order is necessarily soluble is represented by either of the following inclusions
$$ (\ZZ/2\ZZ\longrightarrow 0)^{\rtt l} \subset (\AbKer)^\lr
$$
$$ (2\ZZ\lra \ZZ)^{\rtt r} \cap (0\lra *)^\lr \subset (0\lra \Bbb Q/\ZZ)^\lr
$$
calculated in the category of Finite Groups.
\end{coro}
%
%
%
%

\subsection{$p$-, $p'$-, and $p,p'$-core as an example of a weak factorisation system}
Axiom M2 of a Quillen model category requires that each morphism $A\lra B$ decomposes as 
$$ A\xra{(c)} \cdot \xra{(f)} B$$
where $(c)$ and $(f)$ are orthogonal to each other. 
These decomposition give rise to weak factorisation systems whose existence is proven by the Quillen small object argument.

There are somewhat similar decompositions in group theory. 

That ``each group admits a surjection from a free group'' can be denoted as follows;  each morphism $0\lra G$ admits a decomposition
 $$0\xra{(0\longrightarrow \ZZ)^{rl}} \cdot \xra{(0\longrightarrow \ZZ)^{r}} G $$
in this notation, we think of the Quillen orthogonals as {\em labels} put on arrows, hence 
the notation means that the homomorphisms belong to the corresponding Quillen orthogonals. 

In a finite group, the descending derived series stabilises at a perfect subgroup $P=[P,P]$ (its perfect core) which is characteristic, 
corresponds to the unique decomposition of form 
$$ H \xra{(\AbKer)^{\rtt l}} \cdot \xra{(\AbKer)^{\rtt lr}} G
$$
of a morphism into a map with a perfect kernel $P$, and a map with a soluble kernel. 

Note these decompositions are analogous to decompositions appearing in weak factorisation systems proved by the Quillen small object argument.

\begin{coro}[$p$-core, $p'$-core, $p,p'$-core]\label{coro:2} In the category of Finite Groups, 
\begin{itemize}
\item the $p$-core $O_{p}(G)$ of $G$, i.e.~the largest normal $p$-subgroup of $G$, 
is the group appearing in the unique decomposition of form
  $$ G \xra {(\ZpZ\lra 0)^\rr} G/O_p(G) \xra{(\ZpZ\lra 0)^{\rtt rrl}} 0$$

\item the $p'$-core $O_{p'}(G)$ of $G$, i.e.~the largest normal $p'$-subgroup of $G$, 
is the group appearing in the unique decomposition of form  
   $$ G \xra {(\ZpZ\lra 0)^{\rtt r}} G/O_{p'}(G) \xra{(\ZpZ\lra 0)^{\rtt rl}} 0$$

\item the $p,p'$-core $O_{p,p'}(G)=O_p(G/O_{p'}(G))$ of $G$ 
is the group appearing in the unique decomposition of form  
   $$ G \xra {(\ZpZ\lra 0)^{\rtt r}} G/O_{p'}(G)  \xra {(\ZpZ\lra 0)^{\rtt rr}} G/O_p(G/O_{p'}(G))   \xra{(\ZpZ\lra 0)^{\rtt rl}} 0$$
\end{itemize}
\end{coro}

We end with a couple of test questions suggested by Bob Oliver (private communication).
\begin{que}
Can the Sylow theorems (only existence and uniqueness of Sylow subgroups)
be proven using the characterization of p-groups by Corollary~\ref{coro:1}(4) ?

Could the components of a finite group, and their properties (commute
pairwise, commute with normal p-subgroups) be characterized and proven to exist with
help of our reformulations? 
\end{que}

\subsection{$f$-local groups, localisations and nilpotent groups.\label{sec:locals}}
\footnote{We thank S.O.Ivanov for pointing out the notion of $f$-local groups and the  conjecture of Farjoun 
that a localisation of a nilpotent group is nilpotent [AIP].} 
\def\urtt{\, !\!\rightthreetimes\,}
\def\uurtt{\,!\!\!\!\rightthreetimes\,}
\def\Groups{{\text{Groups}}}
\def\Id{{\text{Id}}}

%

Let $ f \uurtt g$ denote the {\em unique} lifting property. 
For a morphism $f$ of groups, a group $A$ is called $f$-local iff $ f \uurtt A\lra 0$. Under some assumptions, 
each morphisms $H\xra g G$ of groups decomposes as  
$$ H \xra{(f)^{\urtt rl}} \cdot \xra{(f)^{\urtt r}} G$$  
A diagram chasing argument shows that whenever such a decomposition always exists, 
there is a functor  $L=L_f: \Groups\lra \Groups $ defined by
$$H \xra{(f)^{\urtt rl}} L(H) \xra{(f)^{\urtt r}} 0,$$
a natural transformation $\eta: \Id \lra  L: \Groups\lra\Groups$ which induces 
isomorphisms $\eta_G:L\, G \xra{(iso)} LL\, G$, $\eta_{LG}=L(G\xra{\eta_X}LG): L\,G\lra LL\,G$.
A functor with these data is called an {idempotent monad} or a {\em localisation}, and by [CSS] Vopenka principle implies 
that any localisation is of this form. See [AIP] for details and references.  

Our notation allows to express a property closely related to the conjecture 
of Farjoun 
that 
the localisation of a nilpotent group is nilpotent, as follows;
see [AIP] and references therein for a discussion of this conjectures.

Note the diagram has a symmetry: it mentions the diagonal map $H\lra H\times H$.

\begin{conj}[Farjoun] 
The following diagram holds for any property (class) of homomorphisms $L$. 
$$
\xymatrix@C+2pc{ {H} \ar[r]^{(L)^{\urtt rl}} \ar@{->}[d]|{(0\lra *)^{\urtt lr}}  & {H_L} \ar[d]|{\therefore(0\lra *)^{\urtt lr}}  \ar[r]^-{{(L)^{\urtt r}}} & 0 \ar[d] \\
 {H\times H} \ar[r]^{(L)^{\urtt rl}}  & {H_L\times H_L}   \ar[r]^--{{(L)^{\urtt r}}} & 0
}
\ \ \ \ \xymatrix@C+2pc{ {H} \ar[r]^{(L)^{\urtt rl}} \ar@{->}[d]|{(0\lra \Ab)^{\urtt lr}}  & {H_L} \ar[d]|-{\therefore(0\lra \Ab)^{\urtt lr}}  \ar[r]^-{{(L)^{\urtt r}}} & 0 \ar[d] \\
 {H\times H} \ar[r]^{(L)^{\urtt rl}}  & {H_L\times H_L}   \ar[r]^--{{(L)^{\urtt r}}} & 0
}
$$
In the diagram above, ``$\therefore(label)$" reads as:
given a (valid) diagram whose arrows have properties indicate by their labels, the arrow marked by $\therefore$ has the property indicated by its label.
See~Fig.~1 and Corollary~2.3(2) for explanations and details. 
\end{conj}

Our notation suggests the following modifications of the conjecture.

\begin{que} Does it hold for each morphism $H\lra G$ of groups and any 
homomorphism $f$: 
$$\xymatrix@C+5pc{ {H} \ar[r]^{(f)^{\urtt rl}} \ar@{->}[d]|{(0\lra \Ab)^{\urtt lr}}  & {H_f} \ar[d]|-{\therefore(0\lra \Ab)^{\urtt lr}}  \ar[r]^-{{(f)^{\urtt r}}} & 0 \ar[d] \\
 {G} \ar[r]^{(f)^{\urtt rl}}  & {G_f}   \ar[r]^--{{(f)^{\urtt r}}} & 0
}
$$
\end{que}

\begin{que} Does it hold for any diagonal morphism $H\lra H\times H$ of groups, any 
properties (classes) $L$ and $P$ of homomorphisms: 

$\xymatrix@C+3pc{ {H} \ar[r]^{(L)^{\urtt rl}} \ar@{->}[d]|{(P)^{\urtt lr}}  & {H_L} \ar[d]|-{\therefore(P)^{\urtt lr}}  \ar[r]^-{{(L)^{\urtt r}}} & 0 \ar[d] \\
 {H} \ar[r]^{(L)^{\urtt rl}}  & {H\times H_L}   \ar[r]^--{{(L)^{\urtt r}}} & 0
}
$
$\xymatrix@C+3pc{ {H} \ar[r]^{(L)^{rl}} \ar@{->}[d]|{(P)^{lr}}  & {H_L} \ar@{..>}[d]|-{\exists\,(P)^{lr}}  \ar[r]^-{{(L)^{r}}} & 0 \ar[d] \\
 {H} \ar[r]^{(L)^{rl}}  & {H\times H_L}   \ar[r]^--{{(L)^{r}}} & 0
}
$
\end{que}

\begin{que} Under what assumptions on morphism $f:H\lra G$, properties $L$ and $P$ 
of homomorphisms it holds: 

$
\xymatrix@C+3pc{ {H} \ar[r]^{(L)^{\urtt rl}} \ar@{->}[d]_{(P)^{\urtt lr}}  & {H_L} \ar[d]|{\therefore(P)^{\urtt lr}}  \ar[r]^-{{(L)^{\urtt r}}} & 0 \ar[d] \\
 {G} \ar[r]^{(L)^{\urtt rl}}  & {G_L}   \ar[r]^--{{(P)^{\urtt r}}} & 0
}$
$\ \ \ \ \ \ 
\xymatrix@C+3pc{ {H} \ar[r]^{(L)^{\rtt rl}} \ar@{->}[d]_{(P)^{\rtt lr}}  & {H_L} \ar@{..>}[d]|{\exists\,(P)^{\rtt lr}}  \ar[r]^-{{(L)^{\rtt r}}} & 0 \ar[d] \\
 {G} \ar[r]^{(L)^{\rtt rl}}  & {G_L}   \ar[r]^--{{(P)^{\rtt r}}} & 0
}$
\end{que}

In an obvious way the notation suggests a large number of similar questions.
The following is only an example, there is little motivation for this particular choice. 
We use this example as an opportunity to use shortened notation.

\begin{que} Under what assumptions on 
properties $\Delta$, $L$ and $P$ 
of homomorphisms it holds: 
$$
\xymatrix@C+3pc{ {\cdot} \ar[r]^{(L)^{\urtt r..rl}} \ar@{->}[d]_{(P)^{\urtt l..lr}}^{(\Delta)}  & {\cdot} \ar[d]|{\therefore(P)^{\urtt l..lr}}  \ar[r]^-{{(L)^{\urtt r..r}}} & \cdot \ar[d]|{(\Delta)} \\
 {\cdot} \ar[r]|{(L)^{\urtt r..rl}}  & {\cdot}   \ar[r]|--{{(P)^{\urtt r..r}}} & \cdot 
}
\ \ \ \ \ \ 
\xymatrix@C+3pc{ {\cdot} \ar[r]^{(L)^{\rtt r..rl}} \ar@{->}[d]_{(P)^{\rtt l..lr}}^{(\Delta)}  & {\cdot} \ar@{..>}[d]|{\exists\,(P)^{\rtt l..lr}}  \ar[r]^-{{(L)^{\rtt r..r}}} & \cdot \ar[d]|{(\Delta)} \\
 {\cdot} \ar[r]|{(L)^{\rtt r..rl}}  & {\cdot}   \ar[r]|--{{(P)^{\rtt r..r}}} & \cdot 
}$$
In particular, when does it hold for $\Delta=\{H\lra H\times H\,:\,H\text{ a group}\}$ the class of diagonal embeddings?
\end{que}

\subsection{Reformulations with less notation}
In this subsection in a verbose manner  we decipher the notation of Quillen negation of the examples below. 
Fig.~3 represents considerations below as diagrams.

There is no non-trivial homomorphism from a group $F$ to $G$, write $F\nto G$, iff  
$$0\lra F\,\rtt 0\lra G\text{ or equivalently }F\lra 0 \rtt G\lra 0.$$
A group $A$ is {\em Abelian} iff 
\[ \left<a,b\right> \,\lra\,\left<a,b:ab=ba\right>  \rtt\,\, A\lra 0
\]
where $\left<a,b\right> \,\lra\,\left<a,b:ab=ba\right>$ is the  abelianisation morphism sending the free group into the Abelian free group on two generators; 
a group $G$ is {\em perfect}, $G=[G,G]$, iff $G\nto A$ for any Abelian group $A$, i.e. 
\[ \left<a,b\right> \,\lra\,\left<a,b:ab=ba\right>  \rtt\,\, A\lra 0\ \implies\ G\lra 0 \rtt A\lra 0\]
equivalently, for an arbitrary homomorphism $g$, 
\[ \left<a,b\right> \,\lra\,\left<a,b:ab=ba\right>  \rtt\,\,g \ \ \implies\ G\lra 0 \rtt\, \,g\,\]
Yet another reformulation is  that, for each group $S$, 
$$
0\lra G \,\rtt\, [S,S]\lra S .$$
In the category of finite or algebraic groups,
a group $H$ is {\em soluble} iff  $G\nto H$ for each perfect group $G$, 
i.e. 
$$0\lra G\,\rtt 0\lra H\text{ or equivalently }G\lra 0 \rtt H\lra 0.$$
Alternatively, a group $H$ is {\em soluble} iff for every homomorphism $f$ it holds
$$ f \,\rtt\, [G,G]\lra G \text{ for each group }G\  \implies\ f \,\rtt 0\lra\, H 
.$$
A prime number $p$ does not divide the number  elements of a finite group $G$ 
iff  $G$ has no element of order $p$, i.e. no element $x\in G$ such that $x^p=1_G$ yet $x^1\neq 1_G,...,x^{p-1}\neq 1_G$,
equivalently $\ZpZ\nto G$, i.e.
$$0\lra \ZpZ\,\rtt 0\lra G\text{ or equivalently }\ZpZ\lra 0 \rtt G\lra 0.$$
A finite group $G$ is a $p$-group, i.e. the number of its elements is a power of a prime number $p$, iff in the category of finite groups 
$$0\lra \ZpZ \rtt 0\lra H \implies 0\lra H\rtt 0\lra G.$$

A group $H$ is the normal closure of the image of $N$, 
i.e.~no proper normal subgroup of $H$ contains the image of $N$, iff for an arbitrary group $G$
$$N\lra H \,\rtt 0 \lra G.$$
A group $D$ is a subnormal subgroup of a finite group $G$ iff  
$$ N\lra H \,\rtt 0 \lra B \text{ for each group }B\ \implies  N\lra H\,\rtt\,D \lra G
$$
i.e. $D\lra G$ right-lifts wrt any map $N\lra H$ such that $H$ is the normal closure of the image of $N$; the lifting property implies that $D\lra G$ is injective. Recall that $D$ is a subnormal subgroup of a finite group $G$ iff there is a finite series of subgroups
$$D=G_0 \vartriangleleft G_1  \vartriangleleft \ldots  \vartriangleleft G_n =G$$ 
such that $G_i$ is normal in $G_{i+1}$, $i=0,\ldots,n-1$.
This is probably the only claim which requires a proof. First notice that if $D$ is normal in $G$ then 
the lifting property holds.
Given a square corresponding to $N\lra H \,\rtt\, D\lra G$, the preimage of $D$ in $H$ 
is a normal subgroup of $H$ containing the image of $N$, hence the preimage of $D$ contains $H$ and the lifting property holds. 
The lifting property is closed under composition, hence it holds for subnormal subgroups as well.
Now assume $D$ is not subnormal in $G$. As $G$ is finite, there is a minimal subnormal subgroup $D'>D$ of $G$.
By construction 
no proper normal subgroup of $D'$ contains  $D$
but the lifting property $D\lra D' \,\rtt\, D\lra G$ fails.

Finally, a finite group $G$ is nilpotent iff the diagonal group $G$ is subnormal in $G\times G$ [Nilp], i.e.~iff
the diagonal map $G\xra \Delta G\times G$, $g\mapsto (g,g)$
 right-lifts wrt any $N\lra H$ such that $H$ is the normal closure of the image of $N$,  
$$ N\lra H \,\rtt 0 \lra B \text{ for each group }B\ \implies  N\lra H\,\rtt\, G\xra \Delta G\times G
.$$

\section{Speculations. Diagrams commuting up to conjugation}

In this section we make some speculation and remarks on ways to extend our notation.

It is useful to consider diagrams which commute {\em up to conjugation}. 
Inner automorphisms have the following properties which are useful in a diagram chasing computation,
and which in fact characterise inner automorphisms among all automorphisms [Inn, Sch]:
\bi\item
An inner automorphism of a group $G$ extends along any group homomorphism $\iota: G\lra H$, 
i.e.~for any $g\in G$, the inner automorphism $G\lra G, x\mapsto gxg^{-1}$ extends 
to an inner automorphism $H\lra H, x\mapsto \iota(g)x\iota(g^{-1})$
\item An inner automorphism of a group $G$ lifts along any surjective group homomorphism $\iota: G\lra H$, 
i.e.~for any $g\in G$, the inner automorphism $G\lra G, x\mapsto gxg^{-1}$ extends    
to an inner automorphism $H\lra H,\ x\mapsto \iota^{-1}(g)x\iota^{-1}(g^{-1})$ 
\ei

\begin{figure}
\begin{center}
\large
$(a)\  \xymatrix{ A \ar[r]^{i} \ar@{->}[d]_f & X\ar@{-->}[r]^\sigma & X^\sigma \ar[d]|--{\sigma g\sigma^{-1}} \\ B \ar[r]|-{j} \ar@{-->}[urr]|{{\tilde j}}& Y \ar@{-->}[r]^\sigma & Y }$
$(b)\  \xymatrix{ P_1 \ar[r] \ar@{->}[d]|{(0\lra\ZpZ)^{\rttt lr}}  & S_p\ar@{-->}[r]^\sigma & S_p^\sigma \ar[d]|--{\sigma g\sigma^{-1}} 
\\ P_2 \ar[r]  \ar@{-->}[urr]  & G \ar[r]^\sigma & G }$
\end{center}
\caption{\label{fig1}\small
(a) 
The definition of a lifting property $f\rttt g$ {\em up to conjugation}.
\newline
(b) A corollary of Sylow theorem: any $p$-subgroup is contained in the Sylow subgroup $S_p$ up to conjugation.
To see this, take $P_1$ to be the trivial group, and note that $P_2$ in ${(0\lra\ZpZ)^{\rttt lr}}$ means $P_2$ is a $p$-group. 
To see that this property holds for the Sylow subgroup, note that $P_1\lra P_2$ in ${(0\lra\ZpZ)^{\rttt lr}}$ implies 
there is a subgroup subgroup series with $\ZpZ$ quotients connecting $P_1$ and $P_2$, hence 
$Card\, P_2/Card\, P_1$ is a power of $p$, hence $Card\,\Imm\,P_2$ is a power of $p$, hence maps to $S_p$
up to conjugation.
}
\end{figure}

Say that two morphisms $A \xra f B$ and $X \xra g Y$ in a category $C$ are {\em orthogonal up to conjugation} and write $f\rttt g$
iff 
for each $i:A\lra X$ and $j:B\lra Y$
making the square commutative, i.e.~$f\circ j=i\circ g$ there is a diagonal arrow $\tilde j:B\lra X$ 
and an inner automorphism $\sigma:Y\lra Y$ such that 
making the total diagram
$A\xra f B \xra {\tilde j} X^\sigma \xra g Y, A\xra i X \xra \sigma X^\sigma , B\xra j Y\xra \sigma Y	$ commutative, i.e.~$f\circ \tilde j=i$ and $\tilde j\circ g=j$ (see Figure 3a).
Define {\em left/right Quillen negation} or {\em left/right Quillen orthogonal} $P^{\rttt l}$, $P^{\rttt r}$ {\em up to conjugation}
in the obvious way. 

Then the corollary of the Sylow theorem that there is a $p$-subgroup which contains any other $p$-subgroup up to conjugation, 
and such a $p$-subgroup is unique up to conjugation, can be expressed as: 
each morphism  $0\lra G$ decomposes as 
$$
0 \xra {(0\lra\ZpZ)^{\rttt lr}} S_p \xra {(0\lra\ZpZ)^{\rttt lrr}} G,
$$
and such decomposition is unique up to conjugation.

\def\rrtt#1#2#3#4#5#6#7{\xymatrix{ {#1} \ar[r]^{} \ar@{->}[d]_{#2} & {#4} \ar[d]_{#5}^{#6} \\ {#3}  \ar[r] \ar@{-->}[ur]^{}& {#7} }}
\begin{figure}
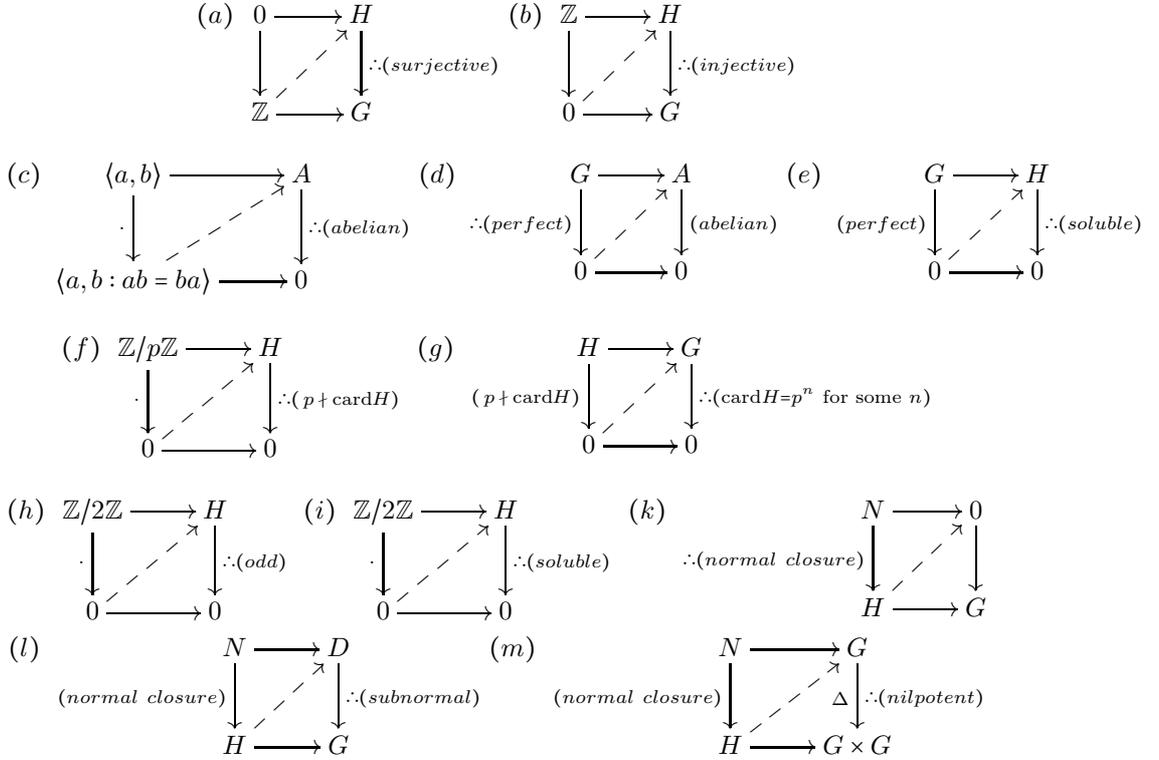

\begin{center}
\small
$(a)\  \rrt  {0} {} {\Bbb Z}   H {\therefore(surjective)} G $
$(b)\  \rrt   {\Bbb Z} {} 0  H {\therefore(injective)} G $
\vskip3mm
$(c)\  \rrt  { \left<a,b\right>} {.} {\left<a,b:ab=ba\right>}   A {\therefore(abelian)} 0 $
$(d)\  \rrt {G} {\therefore(perfect)} {0}  A {(abelian)} 0 $
$(e)\  \rrt G {(perfect)} {0} {H} {\therefore(soluble)} {0}$\ 
\vskip3mm
$(f)\  \rrt {{\Bbb Z}/p{\Bbb Z}} {.} {0} {H} {\therefore(\,p\,\nmid\,\textrm{card}H)} {0}$
$(g)\ \rrt {H} {(\,p\,\nmid\,\textrm{card}H)} {0} {G} {\therefore(\textrm{card}H=p^n\text{ for some }n)} {0}$\\  
\vskip3mm
$(h)\  \rrt {{\Bbb Z}/2{\Bbb Z}} {.} {0} {H} {\therefore(odd)} {0}$
$(i)\ \rrt {{\Bbb Z}/2{\Bbb Z}} {.} {0} {H} {\therefore(soluble)} {0}$
$(k)\ \rrt  {N} {\therefore(normal\ closure)} {H} {0} { } {G}$ 
$(l)\ \rrtt   {N} {(normal\ closure)} {H} {D} {} { \therefore(subnormal)  } {G}$%
 $(m)\ \rrtt   {N} {(normal\ closure)} {H} {G} {\Delta} { \therefore(nilpotent)  } {G\times G}$
\end{center}
\caption{\label{fig5}\footnotesize 
Lifting properties/Quillen negations. Dots $\therefore$ indicate free variables. 
Recall these diagrams represent rules in a diagram chasing calculation
and ``$\therefore(label)$" reads as: 
given a (valid) diagram, add label $(label)$ to the corresponding arrow. 
A diagram is valid
iff for every commutative square of solid arrows
with properties indicated by labels, 
there is a diagonal (dashed) arrow making the total diagram commutative. 
A single dot indicates that the morphism is a constant.\newline
(a) a homomorphism $H\lra G$ is surjective, i.e.~for each $g\in G$ there is $h\in H$ sent to $g$\newline  
(b) a homomorphism $H\lra G$ is injective, i.e.~the kernel of $H\lra G$ is the trivial group\newline
(c) a group is abelian iff each morphism from the free group of two generators 
 factors through its abelianisation ${\Bbb Z}\times {\Bbb Z}$.\newline
 (d)  a group $G$ is perfect, $G=[G,G]$, iff it admits no non-trivial homomorphism to an abelian group\newline
 (e) a finite group is soluble  iff it admits no non-trivial homomorphism from a perfect  group; 
 more generally, this is true in any category of groups with a good enough dimension theory.\newline
(f) by Cauchy's theorem, a prime $p$ divides the number of elements of a finite group $G$ 
iff the group contains an element $e,e^p=1, e\neq1$ of order $p$\newline
(f) a group has order $p^n$ for some $n$ iff iff the group contains no element $e,e^l=1, e\neq1$ of order $l$ prime to $p$\newline
(h)  by Cauchy's theorem, a finite group has an odd number of elements iff it contains no involution $e,e^2=1, e\neq1$\newline
(i) The Feit-Thompson theorem says that each group of odd order is soluble,~i.e.~it says that this diagram chasing 
 rule is valid in the category of finite groups. Note that it is not a definition of the label unlike the other  
 lifting properties.\newline
(k) a group $H$ is the normal closure of the image of $N$ iff $N\lra H \,\rtt 0 \lra G$ for an arbitrary group $G$\newline
(l) $D\lra G$ is injective and the subgroup $D$ is a subnormal subgroup 
of a finite group $G$ iff  $D \lra G$
 right-lifts wrt any map $N\lra H$ such that $H$ is the normal closure of the image of $N$\newline
(m) a group $G$ is nilpotent iff the diagonal map $G\xra \Delta G\times G$, $g\mapsto (g,g)$
 right-lifts wrt any inclusion of a subnormal subgroup $N\lra H$
}\end{figure}

\section*{Acknowledgments and historical remarks.}  
This work is a continuation of [DMG]; early history is given there. 

Examples here were motivated by a discussion with S.Kryzhevich.
I thank Paul Schupp for pointing out the characterisation of inner
automorphisms of [Sch]. I thank M.Bays, K.Pimenov, V.Sosnilo and S.Synchuk for proofreading,
and several students for encouraging and helpful discussions. 

Special thanks are due to M.Bays for helpful discussions.

 I  wish to express my deep thanks to Grigori Mints, to whose memory this paper is dedicated \dots


\begin{thebibliography}{10}

\bibitem[AIP]{AIP}
\newblock Danil Akhtiamov, Sergei O. Ivanov, Fedor Pavutnitskiy.
\newblock Self-derived localizations of groups.
\newblock 2019.
\newblock \url{https://arxiv.org/abs/1905.07612}

%
%
\bibitem[CSS]{CSS}
\newblock C. Casacuberta, D. Scevenels, and J. H. Smith.
\newblock  Implications of large-cardinal principles in homotopical localization.
\newblock Adv. Math. 197 (2005), no. 1, 120–139.

%
%
%


\bibitem[GLP]{GLP}
\newblock Misha Gavrilovich.
\newblock Expressive power of the lifting property in elementary mathematics. A draft, current version.
\newblock \url{http://mishap.sdf.org/mints/expressive-power-of-the-lifting-property.pdf} Arxiv arXiv:1707.06615 (7.17)



\bibitem[DMG]{DMG}
\newblock Misha Gavrilovich,
\newblock Point set topology as diagram chasing computations. Lifting properties as intances of  negation.
\newblock The De Morgan Gazette \ensuremath{5} no.~4 (2014), 23--32,  ISSN 2053-1451
\newblock \url{http://mishap.sdf.org/mints/mints-lifting-property-as-negation-DMG_5_no_4_2014.pdf}

\bibitem[GLZ]{GLZ}
\newblock Misha Gavrilovich, Alexandre Luzgarev, Vladimir Sosnilo.
\newblock A decidable fragment of diagram chasing without automorphisms.
\newblock preprint.
\newblock http://mishap.sdf.org/mints-a-decidable-fragment-of-category-theory-without-automorphisms.pdf


\bibitem[GP]{GP}
\newblock Misha Gavrilovich, Konstantin Pimenov.
\newblock A naive diagram-chasing approach to formalisation of tame topology.
\newblock \url{http://mishap.sdf.org/by:gavrilovich/mintsGE.pdf} 


%
\bibitem[Qui]{Qui}
\newblock D. Quillen. 
\newblock Homotopical Algebra. 
\newblock Lecture Notes in Mathematics, vol. 43. Springer, 1967.


\bibitem[Schupp]{Schupp}
\newblock  P. Schupp. 
\newblock A characterization of inner automorphisms. 
\newblock Proceedings of the American Mathematical Society, vol. 101, n. 2,
     pp. 226-228, 1987.
\newblock  \url{http://www.jstor.org/stable/info/2045986}



\bibitem[Inn]{Inn}
\newblock A wiki on characterisations of inner automorphisms
     of    groups.         Extensible automorphisms conjecture.
\newblock      \url{http://groupprops.subwiki.org/wiki/Extensible automorphisms conjecture}

\bibitem[Nilp]{Nilp}
\newblock  A wiki on characterisations of nilpotent groups. Nilpotent
     groups. 
\newblock \url{http://groupprops.subwiki.org/wiki/Nilpotent group}
\end{thebibliography}
\end{document}